\documentclass[10pt]{amsart}
\usepackage{amsfonts,amssymb,xypic}
\def\qq{\mathbb{Q}}
\def\rr{\mathbb{R}}
\def\zz{\mathbb{Z}}

\def\cc{\mathbb{C}}
\def\pp{\mathbb{P}}
\def\ff{\mathbb{F}}
\def\hh{\mathcal{H}}

\def\ss{\mathcal{S}}
\def\gg{\mathbb{G}}

\def\imtau{\mathrm{Im}\tau}

\def\ii{\mathcal{I}}

\theoremstyle{plain}
\newtheorem{theorem}{Theorem}[section]
\newtheorem{lemma}[theorem]{Lemma}
\newtheorem{proposition}[theorem]{Proposition}

\theoremstyle{remark}
\newtheorem{remark}[theorem]{Remark}
\newtheorem{example}[theorem]{Example}

\theoremstyle{definition}

\newcommand*\isom{%
  \xrightarrow{\sim}%
}

\numberwithin{equation}{section}

\begin{document}

\title{Explicit Mumford isomorphism for hyperelliptic curves}
\author{Robin de Jong}

\begin{abstract}
Using an explicit version of the Mumford isomorphism on the moduli
space of hyperelliptic curves we derive a closed formula for the Arakelov-Green
function of a hyperelliptic Riemann surface evaluated at its Weierstrass points.  
\end{abstract}

\maketitle

\thispagestyle{empty}

\section{Introduction}

The main goal of this paper is to give a formula for the Arakelov-Green function
of a hyperelliptic Riemann surface, evaluated at pairs of Weierstrass points
(cf. Theorem \ref{mainresult} below). This formula generalises a result of Bost in \cite{bo}
dealing with the case that the genus is 2. As an application of our
formula we deduce a symmetric form of a classical identity involving Thetanullwerte and
Jacobian Nullwerte, found originally by Thomae in the 19th century 
(cf. Theorem~\ref{thomae}).

The main idea of our approach is to construct an explicit form of  Mumford's
isomorphism in the case of hyperelliptic curves. We recall that if $p : X \to S$ is a
smooth proper curve with sheaf of relative differentials $\omega$, 
one has a canonical isomorphism $\lambda_1^{\otimes
6n^2+6n+1} \isom \lambda_n$ of invertible sheaves on $S$, ascribed to Mumford
\cite{mu2}; here $n$ is any integer $\geq 1$ and $\lambda_n$ denotes the determinant
sheaf $\det p_* \omega^{\otimes n}$. Later on we will find it more convenient 
to use a different form of Mumford's isomorphism, 
involving Deligne brackets, 
but in order to fix ideas we describe what our results
mean in the present setting. Let $\mu_n$ denote the canonical
trivialising section of $\lambda_n \otimes \lambda_1^{-\otimes 6n^2+6n+1}$
defined by Mumford's isomorphism. In \cite{bsp}, Beilinson and Schechtman
give a formula for $\mu_n$ in the case where $p : X \to S$ is a family of
hyperelliptic curves over the complex numbers. Their result is as follows.
Let $S=\cc^{2g+2} \setminus \{ \textrm{diagonals} \}$ and 
let $p : X \to S$ be the
family of hyperelliptic curves given by
\[ y^2 = f_a(x) = \prod_{i=1}^{2g+2} (x-a_i) \, , \, a = (a_i) \in \cc^{2g+2} \,
, \, a_i \neq a_j \,\, \textrm{if} \,\, i \neq j \, . \] Put $\phi = dx/y \in
H^0(X,\omega)$ and consider the bases $B_n$ of $H^0(X,\omega^{\otimes n})$ given by
\[ B_1 = ( \phi, x\phi, \ldots, x^{g-1}\phi ) \, , \]
\[ B_n = ( \phi^n ,x \phi^n,\ldots,x^{n(g-1)}\phi;
y\phi^n,yx\phi^n,\ldots,yx^{(n-1)(g-1)-2}\phi^n ) \,\, \textrm{for} \,\, n \geq 2
\, . \]
Then we have \[ \mu_n = (\textrm{constant}) \cdot \prod_{(i,j), i \neq j}
(a_i-a_j)^{n(n-1)/2} \cdot \det B_n /(\det B_1)^{\otimes 6n^2+6n+1} \] for $a$ 
running through $S$.
The way we make Mumford's isomorphism explicit is that we are able to calculate the constant
appearing in the above formula for $\mu_n$. In fact it will follow that, 
up to a sign, this constant is equal to
$2^{-(2g+2)n(n-1)}$.

\section{Hyperelliptic curves}

Even though our main result deals with hyperelliptic Riemann surfaces, we need
to consider for the proof hyperelliptic curves over arbitrary base
schemes. Let $g \geq 2$ be an integer, and let $S$ be a scheme. We call a
hyperelliptic curve of genus $g$ over $S$ any smooth, proper curve $p : X \to S$
of genus $g$ which admits an involution $\sigma$ such that for every geometric
point $\overline{s}$ of $S$ the quotient $X_{\overline{s}}/\langle \sigma
\rangle$ is isomorphic to $\pp^1_{k(\overline{s})}$. Once such an involution
exists, it is unique; this is well-known for $S = \mathrm{Spec}(k)$ with $k$ an
algebraically closed field, and follows for the general case by the fact that
$\mathrm{Aut}_S(X)$ is unramified over $S$. If $p : X \to S$ is a hyperelliptic curve, we call $\sigma$ the
hyperelliptic involution of $X/S$. Here are some facts which will be useful
later on.
\begin{proposition} \label{quotient}
The quotient map $X \to X/\langle \sigma \rangle$ is a
finite, faithfully flat $S$-morphism of degree 2 onto a smooth, proper $S$-curve
of genus 0. If $X/\langle \sigma \rangle /S$ admits a section, then $X/\langle
\sigma \rangle$ is $S$-isomorphic to $\pp(V)$ for some locally free sheaf $V$ on
$S$ of rank 2.
\end{proposition}
\begin{proof} See \cite{lons}, Proposition 3.3 and Theorem 5.5.
\end{proof}
Let $\omega$ be the sheaf of relative differentials of $X/S$.
\begin{proposition} \label{canonicalimage}
The image of the canonical morphism $\pi : X \to \pp(p_*
\omega)$ is a smooth, proper $S$-curve of genus 0. Its formation commutes with
base change. Moreover, there exists a closed embedding $j: X/\langle \sigma
\rangle \hookrightarrow \pp(p_*\omega)$ such that $\pi = j \circ h$; here $h$ is
the quotient map $X \to X/\langle \sigma \rangle$. 
\end{proposition}
\begin{proof} See \cite{lons}, Lemma 5.7 and Theorem 5.5.
\end{proof}
The action of $\sigma$ has a fixed point subscheme on $X$, which we denote by
$W$. We call this scheme the Weierstrass subscheme of $X$. It is the closed subscheme defined
locally on an affine open subscheme $U=\mathrm{Spec}(R)$ by the ideal generated
by the set $\{ r - \sigma(r) | r \in R \}$.
\begin{proposition} \label{propWei}
The Weierstrass subscheme $W$ of $X/S$ is the subscheme associated to
a relative Cartier divisor on $X$. It is finite and flat over $S$ of degree
$2g+2$, and its formation commutes with base change. Furthermore, it is \'etale
over a point $s \in S$ if and only if the residue characteristic of $s$ is not
equal to 2.
\end{proposition} 
\begin{proof}
See \cite{lons}, Section 6.
\end{proof}
\begin{example} \label{exam}
Consider the proper, flat genus 2 curve $p:X \to S=\mathrm{Spec}(R)$ with $R =
\zz[1/5]$
given by the affine equation $y^2+x^3y=x$. One may check that it has
good reduction everywhere, and it follows that $p:X \to S$ is a hyperelliptic
curve. Over the ring
$R'=R[\zeta_5,\sqrt[5]{2}]$ it acquires six $\sigma$-invariant
sections with one given by $x=0$ and the others
given by $x=-\zeta_5^k \sqrt[5]{4}$ for $k=1,\ldots,5$. The Weierstrass
subscheme of $X'/R'$ is supported on the images of these sections. It is
clear that they do not meet
over points of residue characteristic $\neq 2$, which verifies that indeed the
Weierstrass subscheme is \'etale over such points. Over a
prime of characteristic 2, 
all $\sigma$-invariant sections meet in the point 
given by $x=0$. The quotient map
$X_{\ff_2} \to X_{\ff_2}/\langle \sigma \rangle \cong
\pp^1_{\ff_2}$ is ramified only in this point.
\end{example}
\begin{remark} In general, if $S$ is the spectrum of a field of
characteristic 2, then the quotient map 
$X \to X/\langle \sigma \rangle $ 
ramifies in at most $g+1$ distinct points.
\end{remark}

\section{A canonical trivialising section of $\lambda_1^{\otimes 8g+4}$}

In this section we study the invertible sheaf $\lambda_1 = \det p_* \omega$ for
a hyperelliptic curve $p : X \to S$.
The following proposition is perhaps well-known.
\begin{proposition} \label{firstlemma}
Suppose that $S$ is a regular integral scheme of generic
characteristic $\neq 2$ and let $p : X \to S$ be a hyperelliptic curve of genus
$g \geq 2$. Then the invertible sheaf $\lambda_1^{\otimes 8g+4}$ has a canonical
trivialising section $\Lambda$. In the case that $S = \mathrm{Spec}(R)$ and that
$X$ has an open subscheme $U =\mathrm{Spec}(E)$ with $E=A[y]/(y^2+ay+b)$, where
$A=R[x]$ and $a,b \in A$, one can write
\[ \Lambda = \left(2^{-(4g+4)} \cdot D \right)^g \cdot \left( \frac{dx}{2y+a}
\wedge \ldots \wedge \frac{x^{g-1} dx}{2y+a} \right)^{\otimes
8g+4} \, , \]
where $D$ is the discriminant in $R$ of the polynomial $a^2-4b$ in $R[x]$.
\end{proposition}
For convenience, we give here the proof; most parts of the argument are taken
from \cite{ka}, Section 6. The statement of Lemma \ref{Kausz3} will be of importance again in
the proof of Proposition~\ref{adjunction}. We start by considering
hyperelliptic curves $p : X \to S$ of genus $g \geq 2$ with $S =
\mathrm{Spec}(R)$ where $R$ is a discrete valuation ring, say with residue
field $k$ and with quotient
field $K$, which we assume to be of characteristic $\neq 2$. 
The canonical quotient map
$R \to k$ will be denoted by $r \mapsto  \bar{r}$.
\begin{lemma} \label{Kausz1} (Cf. \cite{ka}, Lemma 6.1) 
After a finite \'etale surjective base change with a discrete
valuation ring $R'$ dominating $R$, the scheme $X' = X \times_R R'$ can be
covered by open affine
subschemes of the shape $U \cong \mathrm{Spec}(E)$ with $ E =
A[y]/(y^2+ay+b) $, where $A=R'[x]$ and $a,b \in A$, such that the polynomials 
$a^2 -4b$ in $ K'[x]$ are separable of degree $2g+2$ and such that
$\deg a \leq g+1$ and $\deg b \leq 2g+2$. For the reduced
polynomials $\overline{a},\overline{b} \in k'[x]$ one always has $\deg
\overline{a} = g+1$ or $\deg \overline{b} \geq 2g+1$.
\end{lemma}
\begin{proof} Locally in the \'etale topology, any smooth morphism has a
section, and hence by Proposition \ref{quotient}  
after a finite \'etale surjective base change with a
discrete valuation ring $R'$ dominating $R$, one obtains by taking the quotient
under $\sigma$ a finite faithfully flat
$R'$-morphism $h' : X' \to \pp^1_{R'}$ of degree 2. Choose a
point $\infty \in \pp^1_{K'}$ such that $X'_{K'} \to \pp^1_{K'}$
is unramified above $\infty$, and let $x$ be a coordinate on
$V=\pp^1_{K'} - \{ \infty \}$. We can then describe $U = h'^{-1}(V)$ as
$U \cong \mathrm{Spec}(E)$ with $E=A[y]/(y^2+ay+b)$ where
$A=R'[x]$ and $a,b \in A$. If we assume the degree of
$a$ to be minimal, we have $\deg a \leq g+1$ and $\deg b \leq
2g+2$. By Proposition
\ref{propWei}, the Weierstrass subscheme $W$ of $X'/S'$ is
finite and flat over $S'$ of degree $2g+2$. By definition, the
ideal of $W$ is generated by $y-\sigma(y) = 2y+a$ on
$U$. Note that $(2y+a)^2 = a^2-4b$, which defines the norm under
$h'$ of $W$ in $\pp^1_{R'}$. Since this norm is also
finite and flat of degree $2g+2$ over $B'$, and since
$W$ is entirely supported in $U$ by our choice of
$\infty$, we obtain that $\deg(a^2-4b) = 2g+2$. Since the norm of
$W \times_{R'} K'$ in $\pp^1_{K'}$ is \'etale over 
$K'$ by Proposition
\ref{propWei}, the polynomial $a^2-4b$ in $K'[x]$ is separable.
Consider finally the reduced polynomials
$\overline{a},\overline{b} \in k'[x]$. Regarding $y$ as an element
of $k'(X'_{k'})$, we have $\mathrm{div}(y) \geq - \min( \deg
\overline{a}, \frac{1}{2} \deg \overline{b} ) \cdot
h'^*(\overline{\infty})$ by the equation for $y$. On the other
hand it follows from the theorem of Riemann-Roch that $y$ has a pole at both
points of $h'^*(\overline{\infty})$ of order strictly larger than
$g$. This gives the last statement of the lemma.
\end{proof}
\begin{lemma} \label{Kausz2} (Cf. \cite{ka}, Proposition 6.2)
Suppose we have on $X$ an open affine subscheme $U \cong
\mathrm{Spec}(E)$ as in Lemma \ref{Kausz1}. Then the
differentials $x^i dx/(2y+a)$ for $i=0,\ldots,g-1$ are nowhere
vanishing on $U$ and extend to regular global sections of the
sheaf of relative differentials
$\omega$ of $X/S$.
\end{lemma}
\begin{proof} Let $F$ be the polynomial $y^2+a(x)y+b(x) \in A[y]$, and let
$F_x$ and $F_y$ be its derivatives with respect to $x$ and $y$,
respectively. It is readily verified that the morphism
$\Omega_{E/R} = (E dx + E dy)/(F_x dx + F_y dy) \to E$ given by
$dx \mapsto F_y, dy \mapsto -F_x$, is an isomorphism of
$E$-modules. This gives that the differentials $ x^i dx/(2y+a)$
for $i=0,\ldots,g-1$ are nowhere vanishing on $U$. For the second
part of the lemma, it suffices to show that the differentials $x^i
dx/(2y+a)$ for $i=0,\ldots,g-1$ on the generic fiber $U_K$ extend
to global sections of $\Omega^1_{X_K/K}$---but this is well-known
to be true.
\end{proof}
\begin{lemma} \label{Kausz3} (Cf. \cite{ka}, Proposition 6.3)
Suppose we have on $X$ an open affine subscheme $U \cong
\mathrm{Spec}(E)$ as in Lemma \ref{Kausz1}. Let $D$ be the discriminant in $K$
of the polynomial $f = a^2-4b$ in $K[x]$. Then the
modified discriminant $2^{-(4g+4)} \cdot D$ is a
unit of $R$.
\end{lemma}
\begin{proof} In the case that the characteristic of $k$ is $\neq 2$,
this is not hard to see: we know that $W \times_R k$ is \'etale of
degree $2g+2$ by Proposition \ref{propWei}, and hence $f$
remains separable of degree $2g+2$ in $k[x]$ under the reduction
map. So let us assume from now on that the characteristic of $k$
equals 2. If $B$ is any domain, and if $P(T)=\sum_{i=0}^n u_i T^i$
and $Q(T)= \sum_{i=0}^m v_i T^i$ are two polynomials in $B[T]$, we
denote by $R_T^{n,m}(P,Q)$ the resultant in $B$ of $P$ and $Q$.
It satisfies the following property: suppose that at
least one of $u_n,v_m$ is non-zero, and that $B$ is in fact a field. Then
$R_T^{n,m}(P,Q) = 0$ if and only if $P$ and $Q$ have a root in
common in an extension field of $B$. Let $F$ be the polynomial
$y^2+a(x)y + b(x)$ in $A[y]$ with $A=R[x]$, and let $F_x$ and
$F_y$ be its derivatives with respect to $x$ and $y$,
respectively. We set $Q = R_y^{2,1}(F,F_x)$ and $P=
R_y^{2,1}(F,F_y)$ which is $4b-a^2=-f$. Let $H \in R$ be the leading
coefficient of $P$, and put $\Delta = 2^{-(4g+4)} \cdot D$. 
A calculation (for which see for instance
\cite{lock}, Section 1) shows that $ R_x^{2g+2,4g+2}(P,Q) = (H
\cdot \Delta)^2 $. We can read this equation as a formal
identity between certain universal polynomials in the coefficients
of $a(x)$ and $b(x)$. Doing so, we may conclude that $\Delta \in
R$ and that $H^2$ divides $R_x^{2g+2,4g+2}(P,Q)$ in $R$. To show that $\Delta$
is in fact a unit, we distinguish two cases. First we assume that
$\overline{H} \neq 0$. Then $\deg \overline{P} = 2g+2$ and again a
calculation shows that $
R_x^{2g+2,4g+2}(\overline{P},\overline{Q}) = (\overline{H} \cdot
\overline{\Delta} )^2$. The fact that $X_k$ is smooth implies that
$R_x^{2g+2,4g+2}(\overline{P},\overline{Q})$ is non-zero, and
altogether we obtain that $\overline{\Delta}$ is non-zero. Now we
assume that $\overline{H}=0$. Then since $\overline{P} =
\overline{a}^2$ we obtain that $\deg \overline{a} \leq g$ and
hence $\deg \overline{P} \leq 2g$. By Lemma \ref{Kausz1} we have
then $2g+1 \leq \deg \overline{b} \leq 2g+2$. But then from
$2\deg(y)=\deg(\overline{a}y+\overline{b})$ and $\deg(y) > g$, which holds by
the theorem of Riemann-Roch,
it follows that in
fact $\deg \overline{b} = 2g+2$ and hence $\deg \frac{d
\overline{b}}{dx} = 2g$. This
implies that $\deg \overline{Q} = 4g$. A final calculation shows that
$R_x^{2g,4g}(\overline{P},\overline{Q}) = \overline{\Delta}^2$.
Again by smoothness of $X_k$ we may conclude that
$R_x^{2g,4g}(\overline{P},\overline{Q})$ is non-zero. This
finishes the proof.
\end{proof}
\begin{example} Consider once more the curve over $R=\zz[1/5]$
given by the equation $y^2+x^3y=x$, cf. Example \ref{exam} above. 
In the notation from Lemma
\ref{Kausz1} we have $a=x^3, b=-x$. We compute
$D=\mathrm{disc}(x^6+4x)=2^{12}5^5$ so that $\Delta  = 5^5$ which is
indeed a unit in $R$.
\end{example}
\begin{proof}[Proof of Proposition \ref{firstlemma}]
(Cf. \cite{ma}, Proposition 2.7) Again, since locally in the \'etale topology any smooth morphism has a section,
it follows by Proposition \ref{quotient} that after a faithfully flat base change 
the quotient map
$ X \to X /\langle \sigma \rangle$ 
becomes an $S$-morphism onto a $\pp^1_S$. 
Then by Lemma \ref{Kausz1} we may assume that
the scheme $X$ is covered by affine schemes $U \cong
\mathrm{Spec}(E)$ with $E=A[y]/(y^2+ay+b)$ and $A$ a polynomial
ring $R[x]$. For such an affine scheme $U$, consider
$V=\mathrm{Spec}(A)$. In the line bundle $(\det p_*
\omega_{U/V})^{\otimes 8g+4}$ we have a rational section
\[ \Lambda_{U/V} = (2^{-(4g+4)} \cdot D)^g \cdot \left( \frac{dx}{2y+a}
\wedge \ldots \wedge \frac{x^{g-1} dx}{2y+a} \right)^{\otimes
8g+4} \, , \] with $D$ as in Lemma \ref{Kausz3}. One
can check that this section does not depend on any choice of
affine equation $y^2+ay+b$ for $U$, and moreover, these sections
coincide on overlaps. Hence they build a canonical rational
section $\Lambda$ of $\lambda_1^{\otimes
8g+4}$. By Lemma \ref{Kausz2} and Lemma \ref{Kausz3}, this $\Lambda$ 
is a global trivialising section. The general case follows by
faithfully flat descent.
\end{proof}

\section{Adjunction on the Weierstrass subscheme} \label{delignebr}

In this section we recall the formalism of the Deligne bracket \cite{de2}. 
Using this formalism, we construct here a canonical section of a certain
invertible sheaf on the base $S$ of a hyperelliptic curve $p : X \to S$, which
can be seen as a sort of residue map (as in the classical adjunction formula)
for the Weierstrass subscheme of $X/S$. 

Let's start with an arbitrary proper, flat, locally complete intersection curve
$p : X \to S$. Deligne has shown that there exists a natural rule
that associates to any pair $(L,M)$ of invertible sheaves on $X$ an invertible
sheaf $\langle L,M \rangle $ on $S$, such that the following properties are
satisfied: \\
(i) For invertible sheaves $L_1,L_2,M_1,M_2$ on $X$ we 
have canonical isomorphisms $\langle L_1 \otimes
L_2, M \rangle \isom \langle L_1, M \rangle \otimes \langle L_2 ,M
\rangle$ and  
$\langle L, M_1 \otimes M_2 \rangle \isom \langle
L,M_1 \rangle \otimes \langle L,M_2 \rangle$. \\
(ii) For invertible sheaves $L,M$ on $X$ we have a
canonical isomorphism $\langle L,M
\rangle \isom \langle M,L \rangle$. \\
(iii) The formation of the Deligne bracket commutes with
base change, i.e., each cartesian diagram
\[ \xymatrix{ X' \ar[d]^-{p'} \ar[r]^-{u'} & X \ar[d]^-{p} \\
S' \ar[r]^-{u} & S } \]
gives rise to a canonical isomorphism $u^* \langle L, M \rangle \isom
\langle u'^* L, u'^* M \rangle$. \\
(iv) For $P : S \to X$ a section of $p$ and any invertible sheaf $L$ on $X$ 
we have a canonical isomorphism $P^*L \isom \langle O_X(P),L \rangle$. \\
(v) (Adjunction formula) For the sheaf of relative differentials $\omega$ of
$p$ and any section $P: S \to X$ of $p$ we have a canonical
adjunction isomorphism $\langle P,\omega \rangle \isom 
\langle P,P \rangle^{\otimes -1} $. \\
(vi) (Riemann-Roch) Let $L$ be an invertible sheaf on $X$ and let
$\omega$ be the sheaf of relative differentials of $X/S$. Then we have a
canonical isomorphism
\[ \left( \det Rp_* L \right)^{\otimes 2} \isom \langle L,L
\otimes \omega^{\otimes -1} \rangle \otimes (\det Rp_*
\omega)^{\otimes 2} \] of line bundles on $S$, with $\det Rp_*$ denoting
the determinant of cohomology along~$p$. \\
In fact, one can put
\[ \langle L,M \rangle = \det Rp_* (L \otimes M) \otimes (\det Rp_*L)^{-1} \otimes (\det Rp_*
M)^{-1} \otimes (\det Rp_* \omega) \] and then the properties (i)-(vi) can be
checked one by one. Another fact that will be useful later is that if $D$ is a
relative Cartier divisor on $X$ and if $M$ is an invertible sheaf on $X$, one
has a canonical isomorphism
\[ \langle O_X(D), M \rangle \isom \mathrm{Nm}_{D/S}(M|_{D}) \, , \] where
$\mathrm{Nm}_{D/S}$ denotes the norm.

Now let $p : X
\to S$ be a hyperelliptic curve. We will denote here by $W$ the
invertible sheaf associated to the relative Cartier divisor defined by the
Weierstrass subscheme of $X/S$. This change of notation should cause no
confusion. The Deligne bracket that we are interested in is $\langle W, W \otimes \omega
\rangle$ and the statement that we want to prove about it is as follows.
\begin{proposition} \label{adjunction}
Suppose that $S$ is a regular integral scheme of generic
characteristic $\neq 2$ and let $B$ be the branch divisor of $W/S$. Then we have a
canonical isomorphism $\langle W, W \otimes \omega \rangle \isom O_S(B)$.
Furthermore, let $\Xi$ be the rational section of $\langle W,W \otimes \omega
\rangle$ corresponding to the canonical rational section of $O_S(B)$ under this
isomorphism. Then $2^{-(2g+2)} \cdot \Xi$ is a global trivialising section of
$\langle W, W \otimes \omega \rangle$.
\end{proposition}
\begin{proof} By our remarks above, the invertible sheaf $\langle W, W
\otimes \omega \rangle$ is canonically isomorphic to $\mathrm{Nm}((W \otimes
\omega )|_W)$ and this, in turn, is canonically isomorphic to
$\mathrm{Nm}(\omega_{W/S})$ by the adjunction formula. But the
latter is the discriminant of $W/S$, which is
canonically isomorphic to $O_S(B)$, with $B$ the branch divisor of $W/S$. Now
let's look at $2^{-(2g+2)} \cdot \Xi$ as in the statement of the proposition. We
claim that it has neither zeroes nor poles on $S$. First of all remark that it
suffices to place ourselves in the situation where $S=\mathrm{Spec}(R)$ with $R$
a discrete valuation ring whose fraction field $K$ has characteristic $\neq 2$. 
Perhaps after making a faithfully flat cover we can assume that the Weierstrass
subscheme is supported on $2g+2$ sections $W_1,\ldots,W_{2g+2}$ and that the
image of the canonical map $h:X \to X/\langle \sigma \rangle$ is a $\pp^1_R$. 
We assume that the discrete valuation on $R$ is normalised such that 
$v(K^*)=\zz$. 
The valuation $v(\Xi)$ of $\Xi$ at the closed point $s$ of $S$
is then given by the sum $\sum_{k \neq l} (W_k,W_l)$ of the local
intersection multiplicities $(W_k,W_l)$ above $s$ of pairs of
sections $W_k$. Suppose that $W_k$ is given by a
polynomial $x-a_k$, and write $a_k$ as a shorthand for
the corresponding section of $\pp^1_R$. By the projection formula
we have for the local
intersection multiplicities that $4(W_k,W_l)=(2W_k,2W_l)=(h^*
a_k, h^*a_l) = 2(a_k,a_l)$ for each $k \neq l$
hence $(W_k,W_l) = \frac{1}{2} (\alpha_k,\alpha_l)$ for each $k
\neq l$. Now the local intersection multiplicity
$(a_k,a_l)$ above $s$ on $\pp^1_R$  is calculated to be
$v(a_k - a_l)$. This gives that $v(\Xi) =  \sum_{k \neq
l} (W_k,W_l) = \frac{1}{2} \sum_{k \neq l} v(a_k - a_l)
$. By Lemma \ref{Kausz3} we have $ \sum_{k \neq l} v(a_k - a_l)  =
(4g+4)v(2)$ hence  the valuation of 
$ 2^{-(2g+2)} \cdot \Xi$ vanishes at $s$, 
which is what we wanted. The general case follows from this by faithfully flat
descent. 
\end{proof}

\section{Arakelov theory of compact Riemann surfaces} \label{arakelov}

Our main result gives a relation between the Arakelov-Green function of a
hyperelliptic Riemann surface, evaluated at its Weierstrass points, and the
Faltings delta-invariant of that Riemann surface. We introduce these notions in
the present section; for some motivating background and for more 
results we refer to Arakelov's original paper \cite{ar} and Faltings' paper
\cite{fa}. 

We start by fixing a compact Riemann surface $X$ of positive genus $g$. On the space $H^0(X,\omega)$ of holomorphic differential forms we have
a natural hermitian inner product $(\alpha,\beta) \mapsto \frac{i}{2} \int_X
\alpha \wedge
\overline{\beta}$. Let $(\alpha_1,\ldots,\alpha_g)$ be an orthonormal basis for
this inner product. It can be used to build a smooth real (1,1)-form on $X$
given by $\mu = \frac{i}{2g} \sum_{k=1}^g \alpha_k \wedge
\overline{\alpha_k}$. Obviously $\mu$ does not depend on the choice of
orthonormal basis,
and hence is canonical. The Arakelov-Green function of $X$ is now the unique
function $G : X \times X \to \rr_{\geq 0}$ satisfying the following properties
for all $P \in X$:  
\begin{itemize}
\item[(I)] The function $\log G(P,Q)$ is
$C^\infty$ for $Q \neq P$.
\item[(II)] We can write $\log G(P,Q) = \log |z_P(Q)| + f(Q)$
locally about $P$, where $z_P$ is a local coordinate about $P$ and where $f$ is
$C^\infty$ about $P$. 
\item[(III)] We have $\partial_Q
\overline{\partial}_Q
\log G(P,Q)^2 = 2\pi i \mu(Q)$ for $Q \neq P$. 
\item[(IV)] We have $\int_X \log G(P,Q) \mu(Q)=0$. 
\end{itemize}
Existence and uniqueness of $G$ are proved in \cite{ar}. By an application of
Stokes' theorem one finds the symmetry relation $G(P,Q)=G(Q,P)$ for all $P,Q \in
X$. 

An admissible line bundle on $X$ is a line bundle $L$ on $X$ together with a smooth
hermitian metric on $L$ such that the curvature form of $L$ is a multiple of
$\mu$. Using the Arakelov-Green function, one obtains a canonical structure of
admissible line bundle on line bundles of the form $O_X(P)$, with $P$ a point on
$X$, as follows: let $s$ be the tautological section of $O_X(P)$, then
put $\|s\|(Q)=G(P,Q)$ for any $Q \in X$. By property (III) above, the curvature
form of $O_X(P)$ with this metric is equal to $\mu$. Any other admissible metric
on $O_X(P)$ is a constant multiple of the canonical metric; furthermore we get
canonical metrics on line bundles of the form $O_X(D)$ with $D$ a divisor on $X$
by taking tensor products. A very important admissible line bundle is the line
bundle $\omega$ of holomorphic differentials, endowed with its Arakelov metric
$\| \cdot \|_\mathrm{Ar}$; this metric can be defined by insisting that for
every $P$ on $X$, the residue isomorphism $\omega(P)[P] = (\omega \otimes
O_X(P))[P] \isom \cc$ is an isometry, with $\cc$ having its standard euclidean
metric. It is proved in \cite{ar} that this metric is indeed admissible.

For any admissible line bundle $L$ on $X$, Faltings has defined a certain metric
on the determinant of cohomology $\lambda(L) = \det H^0(X,L) \otimes \det
H^1(X,L)^\lor$ of the underlying line bundle (cf. \cite{fa}, Theorem 1). We do
not recall the definition, but mention only that for $L=\omega$, the metric on
$\lambda(L) \cong \det H^0(X,\omega)$ is the one given by the inner product
$(\alpha,\beta) \mapsto \frac{i}{2} \int_X \alpha \wedge \overline{\beta}$
on $H^0(X,\omega)$. It turns out that the Faltings metric on the
determinant of cohomology can be made explicit using theta functions. 
Let $\hh_g$ be the Siegel upper half
space of complex symmetric $g$-by-$g$-matrices with positive
definite imaginary part. Let $\tau \in \hh_g$ 
be a period matrix associated to a symplectic basis of $H_1(X,\zz)$ and
consider the complex torus $J_\tau(X) = \cc^g
/\zz^g + \tau  \zz^g$ associated to $\tau$. 
On $\cc^g$ one has the Riemann theta function
$ \vartheta(z;\tau) = \sum_{n \in \zz^g} \exp(\pi i {}^t n \tau n + 2\pi i {}^t
n z)$, giving rise to an effective divisor $\Theta_0$ and a line bundle 
$O(\Theta_0)$ 
on $J_\tau(X)$. Now consider on the other hand the set
$\mathrm{Pic}_{g-1}(X)$ of divisor classes of degree $g-1$ on $X$. It comes
with a canonical subset $\Theta$ given by the classes of effective divisors. By
the theorem of Abel-Jacobi-Riemann there is a canonical
bijection $u : \mathrm{Pic}_{g-1}(X) \isom J_\tau(X)$ mapping $\Theta$ onto
$\Theta_0$. As a result, we can equip $\mathrm{Pic}_{g-1}(X)$ with the
structure of a compact complex manifold, together with a divisor $\Theta$ and a
line bundle $O(\Theta)$.  

The function $\vartheta$ is not well-defined on
$\mathrm{Pic}_{g-1}(X)$ or $J_\tau(X)$. We can remedy this by putting
$ \|\vartheta\|(z;\tau) = (\det \imtau)^{1/4} \exp(-\pi {}^t y (\imtau)^{-1}
y)|\vartheta(z;\tau)| $, with $y = \mathrm{Im} \, z$. One can check that
$\|\vartheta\|$ descends to a function on $J_\tau(X)$. By our
identification  $\mathrm{Pic}_{g-1}(X) \isom J_\tau(X)$ we obtain
$\|\vartheta\|$ as a function on $\mathrm{Pic}_{g-1}(X)$. It can be checked that
this function is independent of the choice of $\tau$. Note that $\|\vartheta\|$
gives a canonical way to put a metric on the line bundle
$O(\Theta)$ on $\mathrm{Pic}_{g-1}(X)$. 

For any line bundle
$L$ of degree $g-1$ there is a canonical isomorphism $\lambda(L)
\isom O(-\Theta)[L]$, the fiber of $O(-\Theta)$ at the class in
$\mathrm{Pic}_{g-1}(X)$ determined by $L$. Faltings proves in \cite{fa} that
when we give both sides the metrics discussed above, the norm of
this isomorphism is a constant independent of $L$; he writes it as
$\mathrm{e}^{\delta(X)/8}$. The $\delta(X)$ appearing here is the celebrated Faltings
delta-invariant of $X$. An important formula relating $G$ and $\delta$ follows from
these considerations. Again, let $(\alpha_1,\ldots,\alpha_g)$ be an orthonormal
basis of $H^0(X,\omega)$, and let $P_1,\ldots,P_g,Q$ be distinct 
points on $X$. Then the formula
\[ (*) \qquad \| \vartheta \| (P_1+ \cdots + P_g - Q) = 
\mathrm{e}^{-\delta(X)/8} \cdot \frac{ \| \det \alpha_k(P_l)
\|_{\mathrm{Ar}}}{ \prod_{k<l} G(P_k,P_l)} \cdot
 \prod_{k=1}^{g} G(P_k,Q)  \] holds (see \cite{fa}, p.~402). An important
counterpart to this formula was derived by Gu\`ardia \cite{gu}; we will state
a special case of his formula in Section \ref{thomaeidentity} below.

It is possible for $L,M$ 
admissible line bundles on $X$, to endow the invertible sheaves (vector
spaces) $\langle L,M
\rangle$ with natural metrics (called Arakelov metrics here), such that 
all  isomorphisms in (i)-(v)
of Section \ref{delignebr} above become isometries. 
In particular if $L=O_X(P)$ and $M=O_X(Q)$
then $\langle L,M \rangle$ has a certain tautological section $\langle
s_P,s_Q \rangle$ whose norm is just $G(P,Q)$. Faltings' 
metric on the determinant of cohomology has the property that
for all admissible line bundles $L$ and with the canonical Arakelov metrics on
all Deligne brackets of pairs of admissible line bundles, the   
Riemann-Roch isomorphism (vi) is always an isometry.    

\section{Self-intersection of the sheaf of relative differentials}

The purpose of this section is to prove the following proposition.
\begin{proposition} \label{selfint}
Let $p : X \to S$ be a hyperelliptic curve of genus $g \geq
2$ with sheaf of relative differentials $\omega$. If $P,Q$ are $\sigma$-invariant sections of $p$ then
we have a canonical isomorphism 
\[ \langle \omega , \omega \rangle \isom \langle P,Q \rangle^{\otimes -4g(g-1)}
\] of invertible sheaves on $S$, compatible with base change. If 
$B=\mathrm{Spec}(\cc)$, then the above isomorphism is an isometry, provided both
sides are endowed with their canonical Arakelov metrics.
\end{proposition}
We need one lemma, which is a generalisation of Proposition 1 in Section 1.1 of
\cite{bmmb}.
\begin{lemma} \label{oompje}
Let $p:X \to S$ be a hyperelliptic curve of genus $g \geq 2$ with sheaf of
relative differentials $\omega$. For any
$\sigma$-invariant section $P : S \to X$ of $p$ we have a unique
isomorphism \[ \omega \isom O_X( (2g-2)P) \otimes p^* \langle P,P
\rangle^{\otimes -(2g-1)}
\] that induces, by pulling back along $P$, the adjunction
isomorphism $\langle P, \omega \rangle \isom 
\langle P,P \rangle^{\otimes -1}$. The formation of this isomorphism
commutes with base change. If 
$B=\mathrm{Spec}(\cc)$, then the above isomorphism is an isometry, provided both
sides are endowed with their canonical Arakelov metrics.
\end{lemma}
\begin{proof} 
First of all, let $P$ be \emph{any} section of $p$.
Let $h:X \to X / \langle \sigma \rangle$ be the canonical map. We recall
that $X / \langle \sigma \rangle$ is 
a smooth, proper $S$-curve of genus 0. Let $q:X / \langle \sigma \rangle \to
S$ be its structure morphism. By composing $P$ with $h$ we obtain a
section $Q$ of $q$, and as a result we can write 
$X / \langle \sigma \rangle \cong \pp(V)$ for
some locally free sheaf $V$ of rank 2 on $B$. 
On the other hand, consider the
canonical morphism $\pi:X \to \pp(p_* \omega)$. This gives us a
natural isomorphism $\omega \cong \pi^*(O_{\pp(p_*\omega)}(1))$. Let
$j : X/ \langle \sigma \rangle \hookrightarrow \pp(p_* \omega)$ be the 
closed embedding given by
Proposition \ref{canonicalimage}. Passing to a faithfully flat cover, we get that $j$ is
isomorphic to a Veronese embedding $\pp^1 \hookrightarrow \pp^{g-1}$ (cf.
\cite{lons}, Remark 5.11), and hence, 
using a faithfully flat descent argument, one has a natural
isomorphism $j^*(O_{\pp(p_*\omega)}(1)) \cong O_{\pp(V)}(g-1)$. By
well-known properties of projective bundles 
there exists a unique invertible sheaf $L$ on $S$ such
that $O_{\pp(V)}(g-1) \cong O_{\pp(V)}((g-1)\cdot Q) \otimes
q^*L$. By pulling back along $h$, we find a natural isomorphism
$\omega \isom O_X( (g-1) \cdot
(P+\sigma(P))) \otimes p^*L$. In the special case where $P$ is
$\sigma$-invariant, this leads to a natural isomorphism
$\omega \isom O_X((2g-2)P) \otimes p^*L$.
Pulling back along $P$ we find that $L \cong \langle \omega,P
\rangle \otimes \langle P,P \rangle^{\otimes -(2g-2)}$ and with
the adjunction isomorphism $\langle P,P \rangle \cong \langle -P,
\omega \rangle$ then finally $L \cong \langle P,P \rangle^{\otimes
-(2g-1)}$. It is now clear that we have an isomorphism $\omega
\isom O_X( (2g-2)P) \otimes p^* \langle P,P
\rangle^{\otimes -(2g-1)}$ that induces by pulling back along $P$
an isomorphism $\langle P, \omega \rangle \isom \langle P,P \rangle^{\otimes -1}$. Possibly after multiplying with a
unique global section of $O_S^*$, we can establish that the latter isomorphism
be the adjunction isomorphism. The commutativity with base change is
clear from the general base change properties of $\omega$ and of the Deligne bracket. If $B=\mathrm{Spec}(\cc)$ then our isomorphism
multiplies the Arakelov metrics by a constant because both sides are admissible
and hence have the same curvature form. As the adjunction isomorphism is an
isometry, our isomorphism is an isometry at $P$, and hence everywhere.
\end{proof}
The proof of Proposition \ref{selfint} is strongly inspired by the proof of Proposition 2 in
Section 1.2 of \cite{bmmb}.
\begin{proof}[Proof of Proposition \ref{selfint}]
By Lemma \ref{oompje}, we have canonical isomorphisms \[ \omega
\isom O_X( (2g-2)P) \otimes p^* \langle P,P
\rangle^{\otimes -(2g-1)}\] 
and 
\[\omega \isom
O_X( (2g-2)Q) \otimes p^* \langle Q,Q \rangle^{\otimes-(2g-1)} \, . \] 
It follows that $O_X( (2g-2)(P-Q))$ comes from the
base, say $O_X((2g-2)(P-Q)) \isom p^*L$, 
and hence 
\[ \langle (2g-2)(P-Q),P-Q \rangle \isom P^*p^*L \otimes
Q^*p^*L^{\otimes -1} = L \otimes L^{\otimes -1} \] 
is canonically trivial on $S$. Expanding, we get a canonical isomorphism
\[  \langle P,P \rangle^{\otimes 2g-2} \otimes \langle Q,Q
\rangle^{\otimes 2g-2} \isom \langle P,Q
\rangle^{\otimes 2(2g-2)} \] 
of invertible sheaves on $S$. Expanding next
the right hand member of the canonical isomorphism 
\[ \langle \omega, \omega \rangle \isom \langle O_X(
(2g-2)P) \otimes p^* \langle P,P \rangle^{\otimes -(2g-1)}, O_X(
(2g-2)Q) \otimes p^* \langle Q,Q \rangle^{\otimes -(2g-1)}
\rangle \] 
gives the result. The commutativity with base change
is clear. The statement on the norm follows since all the isomorphisms above are
isometries. This is clear from Lemma \ref{oompje}, except possibly for the
isomorphism $ \langle P,P \rangle^{\otimes 2g-2} \otimes \langle Q,Q
\rangle^{\otimes 2g-2} \isom \langle P,Q
\rangle^{\otimes 2(2g-2)} $. But here the statement follows since 
$O_X( (2g-2)(P-Q))$ comes from the base, and hence its Arakelov metric is
constant. By pulling back along $P$ and along $Q$ this constant is cancelled
away, resulting in the trivial metric on $\langle (2g-2)(P-Q),P-Q \rangle$ 
under its canonical trivialisation.
\end{proof}

\section{Explicit Mumford isomorphism}

Let $p: X \to S$ be a smooth, proper curve 
with sheaf of relative differentials
$\omega$. As was mentioned in the Introduction, we have a canonical isomorphism
$\lambda_1^{\otimes 6n^2+6n+1} \isom \lambda_n$ for any integer $n \geq 1$,
where $\lambda_n$ is defined to be the determinant sheaf $\det p_*
\omega^{\otimes n}$. By
Serre duality, this sheaf equals the determinant of cohomology $\det Rp_*
\omega^{\otimes n}$ of $\omega^{\otimes n}$. Taking $n=2$ and applying the
Riemann-Roch isomorphism of Section \ref{delignebr} we obtain a canonical isomorphism
\[ (M) \qquad \mu : \lambda_1^{\otimes 12} \isom \langle \omega, \omega \rangle \, . \]
We have the following result on the norm of $\mu$.
\begin{proposition} (Faltings \cite{fa}, Moret-Bailly \cite{mo1}) 
\label{mumfordnorm}
Assume that $S
= \mathrm{Spec}(\cc)$ and endow both sides of the isomorphism (M) with their
canonical Arakelov metrics. Let $g$ be the genus of $X$. Then the norm of $\mu$ is equal to $(2\pi)^{-4g}
\mathrm{e}^{\delta(X)}$ where $\delta(X)$ is the Faltings delta-invariant of $X$ as in
Section \ref{arakelov}.
\end{proposition}
Now let's consider the case that $p : X \to S$ is a hyperelliptic curve. Using
the results of Section \ref{delignebr} we can identify a certain power of 
$\langle \omega, \omega \rangle$  with a certain power of $\langle W, W \otimes
\omega \rangle$, where $W$ is the
invertible sheaf associated to the Weierstrass subscheme as in Section
\ref{delignebr}. 
Applying the Mumford isomorphism (M), one can thus identify 
a certain power of $\lambda_1$
with a certain power of $\langle W, W \otimes \omega \rangle$.  
The interesting point is
that in this way one can identify a certain power of the
canonical section $\Lambda$, on the one hand, with a certain power of the
canonical section $2^{-(2g+2)} \cdot \Xi$, on the other. More precisely, one has
the following result.
\begin{theorem} \label{canonicalisom}
Let $p : X \to S$ be a hyperelliptic curve of genus $g \geq 2$
with $S$ a regular integral scheme of generic characteristic $\neq 2$
and suppose that there exist $2g+2$ distinct $\sigma$-invariant sections. Then one
has a canonical isomorphism
\[ \lambda_1^{\otimes 12(8g+4)(4g^2+6g+2)} \isom \langle W, W \otimes \omega
\rangle^{\otimes -4g(g-1)(8g+4)} \, . \] This isomorphism maps
$\Lambda^{\otimes 12(4g^2+6g+2)}$ to $(2^{-(2g+2)} \cdot \Xi)^{\otimes
-4g(g-1)(8g+4)}$, up to a sign. 
In the case that $S = \mathrm{Spec}(\cc)$, the isomorphism has norm
$\left((2\pi)^{-4g}\mathrm{e}^{\delta(X)}\right)^{(8g+4)(4g^2+6g+2)}$, if both
sides are equipped with their canonical Arakelov metrics. 
\end{theorem}
\begin{proof} Let $P,Q$ be distinct $\sigma$-invariant 
sections of $X \to S$. By Proposition \ref{selfint} one
has a canonical isomorphism $\langle \omega, \omega \rangle \isom \langle P,Q
\rangle^{\otimes -4g(g-1)}$, which is an isometry for the canonical Arakelov
metrics. Using the adjunction formula for the Deligne bracket one obtains from
this a canonical isomorphism $\langle \omega, \omega \rangle^{\otimes 4g^2+6g+2}
\isom \langle W, W \otimes \omega \rangle^{\otimes -4g(g-1)}$ which is again an
isometry for the Arakelov metrics. Applying the Mumford isomorphism (M) one gets
a canonical isomorphism $\lambda_1^{\otimes 12(4g^2+6g+2)} \isom 
\langle W, W \otimes \omega \rangle^{\otimes -4g(g-1)}$ having norm
$ \left((2\pi)^{-4g}\mathrm{e}^{\delta(X)}\right)^{4g^2+6g+2}$ by Proposition
\ref{mumfordnorm}. The
required isomorphism and the statement on its norm follow from this by raising
to the $(8g+4)$-th power. Now as to the sections on both sides, recall from
Proposition \ref{firstlemma} that $\Lambda$ is a canonical trivialising section of
$\lambda_1^{\otimes 8g+4}$. On the other hand, by Proposition \ref{adjunction} we have that 
$2^{-(2g+2)} \cdot \Xi$ is a
canonical trivialising section of $\langle W, W \otimes \omega \rangle$. The
proof of the theorem is therefore completed by the following proposition.
\end{proof}
\begin{proposition} (Cf. \cite{hr}, Lemma 2.1) \label{H0} 
Let $\ii_g$ be the stack of
hyperelliptic curves of genus $g \geq 2$. Then $H^0(\ii_g ,
\gg_\mathrm{m})=\{-1,+1\}$.
\end{proposition} 
\begin{proof}
We note that we can describe $\ii_g \otimes \cc$ as the space
of $(2g+2)$-tuples of distinct points on $\pp^1$ modulo projective equivalence.
More precisely one has $\ii_g \otimes \cc = ((\pp^1 \setminus \{
0,1,\infty\})^{2g-1} \setminus \{ \mathrm{diagonals} \})/S_{2g+2}$ where
$S_{2g+2}$ is the symmetric group acting by permutation on $2g+2$ points on
$\pp^1$. According to Theorem 10.6 of \cite{hainmcph} the first homology of 
$(\pp^1 \setminus \{
0,1,\infty\})^{2g-1} \setminus \{ \mathrm{diagonals} \} $ is
isomorphic to the irreducible representation of $S_{2g+2}$
corresponding to the partition $\{2g,2\}$ of $2g+2$; in particular
it does not contain a trivial representation of $S_{2g+2}$. This
proves that $H_1(\ii_g \otimes \cc,\qq)$ is trivial, and hence
$H^0(\ii_g \otimes \cc,\gg_\mathrm{m})=\cc^*$. The statement that 
$H^0(\ii_g,
\gg_\textrm{m}) = \{ -1, +1 \} $ follows from this since $\ii_g \to
\mathrm{Spec}(\zz)$ is smooth and surjective.
\end{proof}

\section{Arakelov-Green function at Weierstrass points}

In this section we derive from Theorem \ref{canonicalisom} our main result, which is an
expression for the Arakelov-Green function of a hyperelliptic Riemann surface,
evaluated at its Weierstrass points, in terms of the discriminant of that
surface and its Faltings delta-invariant. Our formula can be seen as a
generalisation of a formula in Proposition 4 of \cite{bo}, which deals with the
special case of Riemann surfaces of genus 2.

Before we state the theorem, we need to introduce the discriminant. Let $g \geq
2$ be an integer and let again $\hh_g$ be the Siegel upper half space.
For vectors $\eta', \eta'' \in \frac{1}{2} \zz^g$ (viewed as column vectors) 
we have on $\cc^g \times
\hh_g$ a theta function $\vartheta[\eta]$ with theta characteristic $\eta =
(\eta',\eta'')$ given by
\[ \vartheta[\eta](z;\tau) = \sum_{n \in \zz^g} \exp( \pi i  {}^t (n+\eta')
\tau  (n+\eta') + 2\pi i  {}^t (n+\eta') (z+\eta'')) \, . \] 
For a given theta characteristic
$\eta$, the corresponding theta function $\vartheta[\eta](z;\tau)$
is either odd or even as a function of $z$. We call the 
theta characteristic $\eta$ odd if the corresponding theta
function $\vartheta[\eta](z;\tau)$ is odd, and even if the
corresponding theta function $\vartheta[\eta](z;\tau)$ is even.

Now let $X$ be a hyperelliptic Riemann surface of genus $g$. We fix an ordering
$W_1,\ldots,W_{2g+2}$ of its Weierstrass points. As is explained in \cite{mu},
Chapter IIIa, this induces a canonical symplectic basis of $H_1(X,\zz)$. Next
choose a coordinate $x$ on $\pp^1$ which puts $W_{2g+2}$ at infinity. This gives
us an affine equation $y^2=f(x)$ of $X$, with $f$ monic ands separable of degree
$2g+1$. Denote by $\mu_1,\ldots,\mu_g$ the holomorphic differentials on $X$ given in coordinates
by $\mu_1=dx/2y,\ldots,\mu_g=x^{g-1}dx/2y$ and denote by $(\mu|\mu')$ the
period matrix of $\mu_1,\ldots,\mu_g$ on the canonical symplectic basis
of homology fixed by our ordering of the Weierstrass points. 
The matrix $\mu$ is invertible and we put $\tau = \mu^{-1} \mu'$.
This matrix lies in $\hh_g$ and we form from it the complex torus
$J_\tau(X)=\cc^g/\zz^g +\tau \zz^g$. Recall from Section
\ref{arakelov} the Abel-Jacobi-Riemann map
$ u: \mathrm{Pic}_{g-1}(X) \isom J_\tau(X) $ identifying the subset
$\Theta$ of classes of effective divisors of degree $g-1$ with the zero locus of
the Riemann theta function $\vartheta(z;\tau)=\sum_{n \in \zz^g} \exp(\pi i
{}^tn \tau n + 2\pi i {}^tnz)$. It is well-known that this map satisfies
$u([K_X-D])=-u([D])$ for all divisors $D$ of degree $g-1$; here $K_X$ denotes a
canonical divisor on $X$. We obtain a bijection
\[ \{ \textrm{classes of $D$ with $2D \sim K_X$} \} \isom
J_\tau(X)[2] \] and hence a bijection 
\[ \{ \textrm{classes of $D$ with $2D \sim K_X$} \} \isom
\{ \textrm{classes mod $\zz^g \times \zz^g$ of theta characteristics} \} \, \]
given by $[D] \mapsto [(\eta',\eta'')]$ if $u([D])=[\eta'+\tau \cdot \eta'']$ on
$J_\tau(X)$. Using the Weierstrass points of $X$, it is easy to produce divisors
$D$ with $2D \sim K_X$ (we call such divisors semi-canonical divisors for
short). Indeed, let $W$ be any Weierstrass point and let $E$ be a divisor from
the hyperelliptic pencil on $X$; then we have $2W \sim E$. But also we have
$(g-1)E \sim K_X$ hence any divisor of degree $g-1$ with support on the
Weierstrass points is semi-canonical. 

We start here by considering 
semi-canonical divisors of the form 
$W_{i_1} + \cdots + W_{i_g} - W_{i_{g+1}} $ for some
subset $\{ i_1,\ldots, i_{g+1} \}$ of cardinality $g+1$ of $\{1,\ldots, 2g+2
\}$. Such divisors have $h^0$ equal to 0, that is, they 
are never linearly equivalent to an effective divisor. 
The remarkable point is that the corresponding theta characteristic depends
only on the set $\{i_1,\ldots,i_{g+1} \}$, and not on $X$. In other words, we
find a canonical map
\begin{align*} \{ & \textrm{subsets $S$ of $\{1,\ldots, 2g+2 \}$ with $\# S = g+1$} \} \\
& \longrightarrow  
\{ \textrm{classes mod $\zz^g \times \zz^g$ of theta characteristics} \} \, . 
\end{align*}
One can prove that this map is 2-to-1; in fact 
$
W_{i_1} + \cdots + W_{i_g} - W_{i_{g+1}} \sim
W_{i'_1} + \cdots + W_{i'_g} - W_{i'_{g+1}} $
if and only if $\{ i_1,\ldots,i_{g+1} \} = \{ i'_1,\ldots, i'_{g+1} \}$ or
$ \{ i_1,\ldots,i_{g+1} \} \cap \{ i'_1,\ldots, i'_{g+1} \} = \varnothing $. Moreover, the theta characteristics in the image are always even. If $S$
is any subset of $\{1,\ldots,2g+2\}$ of cardinality $g+1$, we denote by $\eta_S$
its corresponding theta characteristic. An explicit formula for this
correspondence is given in \cite{mu}, Chapter IIIa, where one finds much more
details on what we have said above. 

Let $\ss$ be the set of subsets of $\{1,\ldots,2g+2\}$ of cardinality $g+1$. We
define on $\hh_g$ the function
\[ \varphi_g(\tau) = \prod_{S \in \ss} \vartheta[\eta_S](0;\tau)^4 \, . \]
According to \cite{lock}, Section 3 the function $\varphi_g(\tau)$  
is a modular form on
$\Gamma_g(2)$
of weight $4r$ where $r={2g+1 \choose g+1 }$. It generalises the usual Jacobi
discriminant modular form in dimension~1. For period matrices $\tau$ which are
associated as above to hyperelliptic Riemann surfaces, the values
$\varphi_g(\tau)$ can be related to the discriminant of a hyperelliptic
equation.
\begin{proposition}
\label{equal} 
Let $X$ be a hyperelliptic Riemann surface of genus $g \geq
2$. Fix an ordering $W_1,\ldots,W_{2g+2}$ of its Weierstrass points. Consider an
equation $y^2=f(x)$ for $X$ with $f$ monic and separable of
degree $2g+1$, putting $W_{2g+2}$ at infinity. 
Let $\mu_k$ for $k=1,\ldots,g$ be the holomorphic differential on $X$ given in
coordinates by $\mu_k = x^{k-1} dx/2y$ and let $(\mu|\mu')$ be the period matrix
of these differentials on the canonical symplectic basis of homology determined
by the chosen ordering of the Weierstrass points. Let $\tau = \mu^{-1}
\mu'$, let $n={2g \choose g+1}$ and let $r={2g+1 \choose g+1}$. Finally let $D$
be the discriminant of $f$. Then  the equality
\[ D^n = \pi^{4gr} (\det \mu)^{-4r} \varphi_g(\tau) \] holds.
\end{proposition}
\begin{proof} See \cite{lock}, Proposition 3.2.
\end{proof}

For a hyperelliptic Riemann surface $X$ of genus $g \geq 2$ we define the
Petersson norm of the modular
discriminant of $X$ to be 
$\|\varphi_g\|(X) = (\det \imtau)^{2r}|\varphi_g(\tau)|$ where
$\tau$ is any period matrix for $X$ formed on a canonical symplectic basis. It
can be checked that the Petersson norm of the modular discriminant of $X$ does not depend on the choice of
this basis, and hence is a (natural and classical) invariant of $X$. It follows from Proposition
\ref{equal} above that it does not vanish. Our main result is now as follows.
\begin{theorem} \label{mainresult}
Let $X$ be a hyperelliptic Riemann surface of genus $g \geq 2$.
Let $m={2g+2 \choose g}$ and 
$n={2g \choose g+1}$. Then we have
\[  \prod_{(W,W')} G(W,W')^{n (g-1)} =
\pi^{-2g(g+2)m} \cdot \mathrm{e}^{-m(g+2)\delta(X)/4}  \cdot \| \varphi_g
\|(X)^{-\frac{3}{2}(g+1)} \, , \] the product running over all
ordered pairs of distinct Weierstrass points of $X$.
\end{theorem} 
\begin{proof} We compute the norms of the sections $\Lambda$ and $\Xi$ for $X$
(considered as a smooth, proper curve over $S=\mathrm{Spec}(\cc)$) and apply the
result of Theorem \ref{canonicalisom}. The formula then drops out. We start with $\Lambda$. As
usual, we fix an ordering $W_1,\ldots,W_{2g+2}$ of the Weierstrass points of $X$ and let $y^2=f(x)$ with
$f$ monic and separable of degree $2g+1$ be an equation for $X$. A small
computation shows that we may write
\[ \Lambda = \left( 2^{-(4g+4)} \cdot D \right)^g \left( \frac{ dx}{y} \wedge \ldots
\wedge \frac{ x^{g-1}dx}{y} \right)^{\otimes 8g+4}  \] for the canonical
trivialising element of $\det H^0(X,\omega)$, 
where $D$ is the discriminant of $f$.
Let $\mu_k$ for $k=1,\ldots,g$ be the holomorphic differential on $X$ given in
coordinates by $\mu_k = x^{k-1} dx/2y$ and let $(\mu|\mu')$ be the period matrix
of these differentials on the canonical symplectic basis of homology determined
by the chosen ordering of the Weierstrass points. Let $\tau = \mu^{-1}
\mu'$, let $r={2g+1 \choose g+1}$ and put
$\Delta_g = 2^{-(4g+4)n} \cdot \varphi_g$. We can then write, by Proposition
\ref{equal}, 
\[ \begin{aligned} \Lambda^{\otimes n} & = \left( 2^{-(4g+4)} \cdot D \right)^{gn} \left( \frac{
dx}{y} \wedge \ldots \wedge \frac{ x^{g-1}dx}{y}
\right)^{\otimes (8g+4)n} \\
& = 2^{-(4g+4)gn} \pi^{4g^2r} (\det \mu)^{-4gr}
\varphi_g(\tau)^g \left( \frac{ dx}{y} \wedge \ldots \wedge
\frac{ x^{g-1}dx}{y} \right)^{\otimes (8g+4)n} \\ 
& = (2\pi)^{4g^2r} (\det \mu)^{-4gr}
\Delta_g(\tau)^g \left( \frac{ dx}{2y} \wedge \ldots \wedge
\frac{ x^{g-1}dx}{2y} \right)^{\otimes (8g+4)n} \, . \end{aligned} \]
Let $J_\tau(X)=\cc^g/\zz^g + \tau \zz^g$, and let $j \colon \det H^0(X,\omega) \isom \det H^0(J_\tau(X), \omega)$ 
be the canonical isomorphism. Letting $z_1,\ldots,z_g$ be the
standard euclidean coordinates on $J_\tau(X)$ we obtain from the
above calculation
\[ j^{\otimes (8g+4)n} (\Lambda^{\otimes n}) = (2\pi)^{4g^2r} \Delta_g(\tau)^g (dz_1 \wedge \ldots \wedge
dz_g)^{\otimes (8g+4)n } \, . \]
It follows that the norm of $\Lambda$ satisfies
\[  \| \Lambda \|^n = (2\pi)^{4g^2r} \|\Delta_g\|(X)^g \, , \] where
$\|\Delta_g\|(X)=2^{-(4g+4)n} \cdot \|\varphi_g\|(X)$; indeed, by definition the
norm of $dz_1 \wedge \ldots \wedge dz_g $ is $\| dz_1 \wedge \ldots \wedge dz_g \| = \sqrt{\det
\imtau}$. Now we consider the section $\Xi$. It has norm
\[ \| \Xi \| = \prod_{(W,W')} G(W,W') \] with the product running over all
ordered pairs of distinct Weierstrass points of $X$. Applying Theorem
\ref{canonicalisom} we
have 
\[  \left( (2\pi)^{-4g}
\mathrm{e}^{\delta(X)} \right)^{(8g+4)(4g^2+6g+2)} \cdot \| \Lambda
\|^{12(4g^2+6g+2)} =  \| 2^{-(2g+2)} \cdot \Xi \|^{-4g(g-1)(8g+4)} \, .
\]
Plugging in the formulas for $\|\Lambda\|$ and $\|\Xi\|$ that we just gave one
obtains the required formula. 
\end{proof}
\begin{remark} In \cite{jong} we constructed two natural invariants $S(X)$ and
$T(X)$ of compact Riemann surfaces $X$, related to the delta-invariant by
the formula $\mathrm{e}^{\delta(X)/4}=S(X)^{-(g-1)/g^2} \cdot T(X)$. Putting $G' =
S(X)^{-1/g^3} \cdot G$ the formula in Theorem \ref{mainresult} can be rewritten as 
\[ \prod_{(W,W')} G'(W,W')^{n (g-1)} =
\pi^{-2g(g+2)m} \cdot T(X)^{-(g+2)m}  \cdot \| \varphi_g
\|(X)^{-\frac{3}{2}(g+1)} \, . \]
In this form our formula is instrumental in the paper \cite{jong2}, where 
a closed formula is given for the
delta-invariant of $X$.
\end{remark}

\section{A classical identity of Thomae} \label{thomaeidentity}

In this final section we combine our Theorem \ref{mainresult} with a formula due to Gu\`ardia
in order to
obtain a symmetric version of an identity found in the 19th century by
Thomae \cite{th}. This identity relates a certain Jacobian
Nullwert to a certain product of Thetanullwerte in the context of hyperelliptic
period matrices. The classical proof of
Thomae's identity can perhaps best be learnt from 
the paper \cite{frob} by Frobenius. Interestingly, in this classical 
proof the heat equation
for the theta function plays a fundamental role. In our approach the heat
equation is circumvented, which perhaps leads to a better 
`algebraic' understanding of Thomae's identity. We remark that the
relations between Jacobian Nullwerte and Thetanullwerte have been studied
extensively by Igusa, see for instance \cite{ig1} and \cite{ig2}, and recently
again by Gu\`ardia in his paper \cite{gu2}.

Let $g \geq 2$ be an integer. 
Let $\eta_1,\ldots,\eta_g$ be $g$ odd theta characteristics in dimension $g$. 
We recall that the
Jacobian Nullwert $J(\eta_1,\ldots,\eta_g)$ in $\eta_1,\ldots,\eta_g$ is defined
to be the
jacobian
\[ J(\eta_1,\ldots,\eta_g)(\tau) = \frac{ \partial(
\vartheta[\eta_1],\ldots,\vartheta[\eta_g]) }{\partial(z_1,\ldots,z_g) }(0;\tau)
\, , \] viewed as a function on $\hh_g$, the Siegel
upper half space. We want to study the values of Jacobian
Nullwerte for period matrices coming from hyperelliptic Riemann surfaces. 
So let $X$ be a hyperelliptic Riemann
surface of genus $g$ and let $\tau$ be a period matrix associated to a
canonical symplectic basis of $X$, given by a certain ordering
$W_1,\ldots,W_{2g+2}$ of its Weierstrass points. 
We recall from Section \ref{arakelov} that in this set-up,
the Abel-Jacobi-Riemann map $u$ induces a canonical bijection
\[ \{ \textrm{classes of semi-canonical divisors} \} \isom
\{ \textrm{classes mod $\zz^g \times \zz^g$ of theta characteristics} \} \, \]
given by $[D] \mapsto [(\eta',\eta'')]$ if $u([D])=[\eta'+\tau \cdot \eta'']$ on
$J_\tau(X)$. Here we want to consider semi-canonical divisors of the form
$W_{i_1} + \cdots + W_{i_{g-1}}$ for subsets $\{i_1,\ldots,i_{g-1}\}$ of
$\{1,\ldots,2g+2\}$ of cardinality $g-1$. Such divisors have $h^0$ equal to 1. 
Again, the remarkable point is that the theta characteristic corresponding to 
$W_{i_1} + \cdots + W_{i_{g-1}}$ depends only on the set
$\{i_1,\ldots,i_{g-1}\}$, and not on $X$. We end up with a canonical map
\begin{align*} \{ & \textrm{subsets $S$ of $\{1,\ldots, 2g+2 \}$ with $\# S = g-1$} \} \\
& \longrightarrow  
\{ \textrm{classes mod $\zz^g \times \zz^g$ of theta characteristics} \} \, . 
\end{align*}
One can prove that this map is 1-to-1, and that the theta characteristics in
the image are always odd. Again, the correspondence can be made explicit; see
again \cite{mu}, Chapter IIIa for the details. 
Now choose a subset $\{i_1,\ldots,i_g\}$ of
$\{1,\ldots,2g+2\}$ of cardinality $g$, and for $k=1,\ldots,g$ let $\eta_k$ be
the odd theta characteristic corresponding to $\{ i_1,\ldots,\widehat{i_k},\ldots,
i_g \}$ by the above canonical map. We put
\[ \|J\|(W_{i_1},\ldots,W_{i_g}) = (\det
\imtau)^{(g+2)/4}|J(\eta_1,\ldots,\eta_g)(\tau)| \, . \]
It can be checked that this only depends on the set $\{W_{i_1},\ldots,W_{i_g}\}$
and not on the chosen ordering of the Weierstrass points. 
We have the following theorem.
\begin{theorem} (Thomae's identity)  \label{thomae} Let $X$ be a hyperelliptic
Riemann surface of genus $g \geq 2$ with Weierstrass points
$W_1,\ldots,W_{2g+2}$. Let $m = {2g+2 \choose g}$. Then we have
\[ \prod_{ \{ i_1, \ldots, i_g \}  } \|J\|(W_{i_1},
\ldots, W_{i_g}) = \pi^{gm} \|\varphi_g\|(X)^{(g+1)/4} \, , \]
where the product runs over the subsets of $\{1,\ldots,2g+2\}$ of
cardinality $g$.
\end{theorem} 
Our proof is basically a combination of Theorem \ref{canonicalisom} with the following 
proposition, which is a special case of the main theorem of 
\cite{gu}. The formula can be obtained from Faltings'
formula (*) by a limiting process, using Riemann's singularity theorem.
\begin{proposition} (Gu\`ardia \cite{gu}) \label{guardia} Let $W_{i_1},\ldots,W_{i_g},W$ be
distinct Weierstrass points of $X$. Then the formula
\[ \|\vartheta\|(W_{i_1}+\cdots+W_{i_g}-W)^{g-1} = \mathrm{e}^{\delta(X)/8}
\cdot \|J\|(W_{i_1},\ldots,W_{i_g}) \cdot \frac{ \prod_{k=1}^g
G(W_{i_k},W)^{g-1} }{ \prod_{k<l} G(W_{i_k},W_{i_l})} \] holds.
\end{proposition}
\begin{proof}[Proof of Theorem \ref{thomae}] We start by taking a set $\{i_1,\ldots,i_g\}$
and taking the product over all $W$ not in $\{W_{i_1},\ldots,W_{i_g} \}$ in
the formula from Proposition \ref{guardia}. This gives
\[ \begin{aligned}
&\prod_{W \notin \{ W_{i_1},\ldots, W_{i_g} \} } \prod_{k=1}^g
G(W_{i_k},W)^{2g-2} \\ &= \mathrm{e}^{-(g+2)\delta(X)/4} \cdot \frac{
\prod_{W \notin \{ W_{i_1},\ldots, W_{i_g} \} } \| \vartheta
\|(W_{i_1}+\cdots+W_{i_g}-W)^{2g-2}}{\| J
\|(W_{i_1},\ldots,W_{i_g})^{2g+4}} \cdot \prod_{k \neq l}
G(W_{i_k},W_{i_l})^{g+2} \, . \end{aligned}
\] Taking the product over all sets  $\{ i_1, \ldots,
i_g \}$ of cardinality $g$ we find
\[ \begin{aligned}
\prod_{(W,W')} & G(W,W')^{n (g-1)}  \\ &= \mathrm{e}^{-m(g+2)\delta(X)/4} 
\cdot \prod_{ \{ i_1, \ldots, i_g \}  } \frac{ \prod_{
W \notin \{W_{i_1}, \ldots, W_{i_g} \} } \| \vartheta \| (W_{i_1}+
\cdots + W_{i_g} - W)^{2g-2} }{ \|J\|(W_{i_1}, \ldots,
W_{i_g})^{2g+4} } \, .\end{aligned} \] {}From our definition of $\|\varphi_g\|(X)$
it follows that
\[ \|\varphi_g\|(X)=
\prod_{ \{ i_1,\ldots,i_{g+1} \} }
\|\vartheta\|(W_{i_1}+\cdots+W_{i_g}-W_{i_{g+1}})^4 \, , \] where
the product runs over the set of subsets of
$\{1,2,\ldots,2g+2\}$ of cardinality $g+1$. This gives
\[ \prod_{ \{ i_1, \ldots, i_g \}  }   \prod_{ W
\notin \{W_{i_1}, \ldots, W_{i_g} \} } \| \vartheta \| (W_{i_1}+
\cdots + W_{i_g} - W)^{2g-2} = \| \varphi_g \|(X)^{(g^2-1)/2} \, .
\] Plugging this in in our previous formula gives
\[  \begin{aligned} \prod_{(W,W')} &G(W,W')^{n (g-1)} = \\
&= \mathrm{e}^{-m(g+2)\delta(X)/4}  \cdot \|\varphi_g\|(X)^{(g^2-1)/2} \cdot
\prod_{ \{ i_1, \ldots, i_g \} } \|J\|(W_{i_1}, \ldots,
W_{i_g})^{-(2g+4)}  \, . \end{aligned}
\]
Comparing this formula with the one in Theorem \ref{mainresult} gives the required
formula.
\end{proof}
It is possible to derive from Theorem \ref{thomae} a statement involving holomorphic
functions on the domain of hyperelliptic period matrices in $\hh_g$. We call a
set $\{\eta_1,\ldots,\eta_g \}$ of odd theta characteristics special if it can
be obtained from a subset of $\{1,\ldots,2g+2\}$ of cardinality $g$ in the way
that we described above. Let $H$ denote the set of special sets of odd theta
characteristics, and let as before $\ss$ denote the set of subsets of
$\{1,\ldots,2g+2\}$ of cardinality $g+1$. 
Then one can deduce from our result that
for period
matrices $\tau$ associated to canonical symplectic bases of hyperelliptic
Riemann surfaces of genus $g$ one has
\[ \prod_{ \{\eta_1,\ldots,\eta_g\} \in H} J(\eta_1,\ldots,\eta_g)(\tau) = \pm
\pi^{gm} \prod_{S \in \ss} \vartheta[\eta_S](0;\tau)^{g+1} \, . \]
Indeed, one observes first that
by dividing left and right of the formula in Theorem \ref{thomae} by an appropriate power
of $\det \imtau$ one gets
\[ \prod_{\{\eta_1,\ldots,\eta_g\} \in H} |J(\eta_1,\ldots,\eta_g)(\tau)| =
\pi^{gm} | \varphi_g(\tau)|^{(g+1)/4} \, . \]
The maximum principle for holomorphic functions allows us then to write
\[ \prod_{\{\eta_1,\ldots,\eta_g\} \in H} J(\eta_1,\ldots,\eta_g)(\tau) =
\varepsilon \pi^{gm} \prod_{S \in \ss} \vartheta[\eta_S](0;\tau)^{g+1} \, , \] where $\varepsilon$ is a complex number of modulus 1
depending only on $g$. Considering the Fourier expansions on left and right as
in \cite{ig1}, pp. 86-88 one finds the value $\varepsilon = \pm 1$.

\subsection*{Acknowledgements} The author wishes to thank  
Riccardo Salvati Manni, Christophe Soul\'e
and Gerard van der Geer for their encouragement and helpful
remarks. He also thanks the Institut des Hautes \'Etudes Scientifiques in
Bures-sur-Yvette, where
a preliminary version of this article was written, for its hospitality during a
visit in October and November 2004.

\vspace{0.5cm}

\noindent Address of the author: \\ \\
Mathematical Institute \\
University of Leiden \\
PO Box 9512 \\
2300 RA Leiden \\
The Netherlands \\  
E-mail:  \verb+rdejong@math.leidenuniv.nl+

\end{document}